\documentclass[journal]{IEEEtran}

\usepackage{url}
\usepackage{float}
\usepackage{hyperref}
\newcommand{\seqnum}[1]{\href{https://oeis.org/#1}{\rm \underline{#1}}}
\usepackage{caption} 
\usepackage{amssymb}
\usepackage{amsmath}
\usepackage{amsthm}
\usepackage{amsfonts}
\usepackage{amscd}
\usepackage{tikz}

\theoremstyle{definition}
\newtheorem{theorem}{Theorem}

\newtheorem{lemma}[theorem]{Lemma}

\begin{document}

\title{On The Number of Inequivalent Monotone Boolean Functions of 9 Variables}

\author{Bart\l{}omiej Pawelski~
\IEEEcompsocitemizethanks{\IEEEcompsocthanksitem B. Pawelski was with the Institute of Informatics, University of Gdansk.}
\thanks{}}

\markboth{}%
{Pawelski \MakeLowercase{\textit{et al.}}: On The Number of Inequivalent Monotone Boolean Function of 9 Variables}

\IEEEtitleabstractindextext{
\begin{abstract}
We provide the first-ever calculation of the number of inequivalent monotone Boolean functions of 9 variables, which is equal to 789,204,635,842,035,040,527,740,846,300,252,680.
\end{abstract}

\begin{IEEEkeywords}
Boolean function, monotone Boolean function, Dedekind number,  inequivalent monotone Boolean functions
\end{IEEEkeywords}}

\maketitle

\IEEEpeerreviewmaketitle
\IEEEdisplaynontitleabstractindextext

\section{Introduction}

\IEEEPARstart{T}{his} paper is a continuation of our previous work \cite{pawelski} on counting inequivalent monotone Boolean functions, in which we provided the first-ever calculation of inequivalent monotone Boolean functions of 8 variables. Our objective in this paper is to calculate the number of inequivalent monotone Boolean functions of 9 variables.

Let $B$ denote the set of two bits \{0,1\} and let $B^n$ denote the set of $n$--element sequences of $B$. A Boolean function is any function $f: B^n \to B$. We have a partial order in $B$: $0 \leq 0$, $0 \leq 1$, and $1 \leq 1$. This partial order induces a partial order on $B^n$: for any two elements $x= (x_1,\dots,x_n)$ and  $ y=(y_1,\dots,y_n)\in B^n$, $x\le y$ if and only if $x_i \leq y_i$ for all $i$.

A Boolean function is said to be \textit{monotone} if for any pair $x, y \in B^n$, when $x \leq y$, it follows that $f(x) \leq f(y)$.

We define $D_n$ as the set of all monotone Boolean functions of $n$ variables. Let $d_n$ represent the cardinality of $D_n$, also known as the $n$--th Dedekind number. Dedekind numbers are described by the \textit{On-Line Encyclopedia of Integer Sequences} (OEIS) sequence \seqnum{A000372} (see Table \ref{tab:valuesd9}).

Two monotone Boolean functions are said to be \textit{equivalent} if the first function can be obtained from the second function through any permutation of input variables. Let $R_n$ represent the set of all equivalence classes of $D_n$, and let $r_n$ denote the cardinality of this set. The values of $r_n$ are described by the OEIS sequence \seqnum{A003182} (see Table \ref{tab:valuesr9}).

In 1985 and 1986, Liu and Hu \cite{liu1, liu2}  calculated $r_n$ for values of $n$ up to 7. Independently, $r_7$ was computed by Stephen and Yusun \cite{SY}. In 2021, the author calculated $r_8$ \cite{pawelski}, and this result was independently reported in 2022 by Cari\'c and \v{Z}ivkovi\'c \cite{cz}.

A significant contribution to the topic is the recent paper by Szepietowski \cite{szepietowski}. Another important and recent result is the calculation of $d_9$, which was independently achieved in 2023 by Jäkel \cite{jakel} and Van Hirtum, De Causmaecker, Goemaere, Kenter, Riebler, Lass, and Plessl \cite{hirtum}.

Our work involves counting the fixed points in $D_9$ under permutations of the 29 cycle types of its input variables (see Table \ref{tab:r9results}). Then, we use Burnside's lemma to calculate $r_9 = 789204635842035040527740846300252680$.

\begin{table}[h] \centering {\renewcommand{\arraystretch}{1.4}
\begin{tabular}{l|l} $n$ & $d_n$ \\ \hline 
0&  2      \\  
1&  3      \\
2&  6       \\  
3&  20       \\  
4&  168       \\
5&  7581       \\  
6&  7828354     \\  
7&  2414682040998  \\
8&  56130437228687557907788 \\
9&  286386577668298411128469151667598498812366 \\ \hline \end{tabular}} \caption{Known values
of $d_n$.} \label{tab:valuesd9} \end{table}

\begin{table}[h] \centering {\renewcommand{\arraystretch}{1.4}
\begin{tabular}{l|l} $n$ & $r_n$ \\ \hline 
0&  2\\  
1&  3\\
2&  5\\  
3&  10\\  
4&  30\\
5&  210\\  
6&  16353\\  
7&  490013148  \\
8&  1392195548889993358\\
9&  789204635842035040527740846300252680 \\ \hline \end{tabular}} \caption{Known values of $r_n$.} \label{tab:valuesr9} \end{table}

\section{Preeliminaries}

\subsection{Burnside's lemma}
Let $S_n$ represent the set of all permutations of the set \{0,1,\ldots,$n$\}. In this context, we treat $S_n$ as the set of all permutations of the $n$ variables within a Boolean function. Each permutation $\pi \in S_n$ also acts upon $B^n$ and $D_n$. We define an element $x$ as a \textit{fixed point} of $\pi$ if it remains unchanged under the action of $\pi$. The set $\Phi_n(\pi)$ contains all fixed points in $D_n$ associated with the permutation $\pi$. 

\begin{table}[H]
\centering
{\renewcommand{\arraystretch}{1.2}
\begin{tabular}{c|c|c|c|c|c|c|c|c}
& 0 & 1 & 2 & 3 & 4 & 5 & 6 & 7 \\
\hline
$(1)$ & $\varnothing$ & $001$ & $010$ & $011$ & $100$ & $101$ & $110$ & $111$  \\
\hline
$(1\;2)$ & $\varnothing$ & $010$ & $001$ & $011$ & $100$ & $110$ & $101$ & $111$  \\
\end{tabular}}
\caption{$B^3$ under $\pi = (12)$.} \label{tab:rg}
\end{table}

\def\drawconnectionsa{\draw (0,0) -- (-2,2);
    \draw (0,0) -- (0,2);
    \draw (0,0) -- (2,2);
    \draw (-2,2) -- (-2,4);
    \draw (-2,2) -- (0,4);
    \draw (0,2) -- (-2,4);
    \draw (0,2) -- (2,4);
    \draw (2,2) -- (0,4);
    \draw (2,2) -- (2,4);
    \draw (-2,4) -- (0,6);
    \draw (0,4) -- (0,6);
    \draw (2,4) -- (0,6);}
\begin{center}
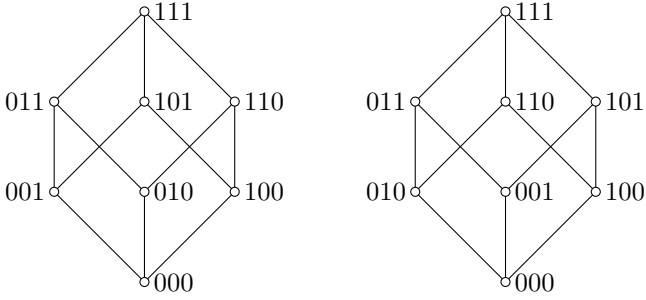


\begin{tikzpicture}[scale=0.6]
\drawconnectionsa{}
\draw[black, fill=white] (0,0) circle (.1cm) node[right] {$000$};
\draw[black, fill=white] (-2,2) circle (.1cm) node[left] {$001$};
\draw[black, fill=white] (0,2) circle (.1cm) node[right] {$010$};
\draw[black, fill=white] (2,2) circle (.1cm) node[right] {$100$};
\draw[black, fill=white] (-2,4) circle (.1cm) node[left] {$011$};
\draw[black, fill=white] (0,4) circle (.1cm) node[right] {$101$};
\draw[black, fill=white] (2,4) circle (.1cm) node[right] {$110$};
\draw[black, fill=white] (0,6) circle (.1cm) node[right] {$111$};

\begin{scope}[shift={(8cm,0)}]
\drawconnectionsa{}
\draw[black, fill=white] (0,0) circle (.1cm) node[right] {$000$};
\draw[black, fill=white] (-2,2) circle (.1cm) node[left] {$010$};
\draw[black, fill=white] (0,2) circle (.1cm) node[right] {$001$};
\draw[black, fill=white] (2,2) circle (.1cm) node[right] {$100$};
\draw[black, fill=white] (-2,4) circle (.1cm) node[left] {$011$};
\draw[black, fill=white] (0,4) circle (.1cm) node[right] {$110$};
\draw[black, fill=white] (2,4) circle (.1cm) node[right] {$101$};
\draw[black, fill=white] (0,6) circle (.1cm) node[right] {$111$};
\end{scope}
\end{tikzpicture}
\captionof{figure}{$B^3$ and $B^3$ under $\pi = (12)$ on Hasse diagrams}
\end{center}

Each permutation $\pi \in S_n$ can be represented as a product of disjoint cycles. Define the
{\it cycle type\/} of $\pi$ to be the tuple of lengths of its disjoint cycles
in increasing order. For example, the cycle type of permutation $\pi =(12)(3456)$ is $(2,4)$, and its total length is 6 \cite{pawelski}.

Using the following application of Burnside's lemma, the $r_n$ can be calculated \cite{cz, pawelski, szepietowski}:

\begin{equation} \label{burnside}
r_n = \frac{1}{n!} \sum_{i=1}^{k} \mu_i \phi(\pi_i),
\end{equation}

where:
\begin{itemize}
\itemsep-0.01em 
  \item[--] $\phi_n(\pi) = |\Phi_n(\pi)|$
  \item[--] $k$ = number of different cycle types in $S_n$
  \item[--] $i$ = index of cycle type
  \item[--] $\mu_i$ = number of $\pi \in S_n$ with cycle type $i$
  \item[--] $\pi_i$ = representative $\pi \in S_n$ with cycle type $i$
\end{itemize}

The first application of Burnside's lemma for calculating $r_n$ that we found in the literature was by Liu and Hu \cite{liu1, liu2}, who used it in 1985 and 1986 to determine $r_n$ for values of $n$ up to 7.
Computing $\phi_9(\pi)$ for all $\pi \in S_9$ (except identity) requires significantly less computational effort compared to calculating $d_9$, making it calculable using resources at our hand. 

\subsection{Posets}

A binary relation that is reflexive, antisymmetric, and transitive, when combined with a set $P$, forms a \textit{partially ordered set}, or simply a \textit{poset}. In the introduction, we defined  partial orders in $B$ and $B^n$  making them both posets.

An incidence matrix of a poset is a binary matrix that represents the partial order relation between elements in the poset. For a poset $(X, \le)$, its incidence matrix $M$ has rows indexed by elements in $X$ and columns indexed by elements in $X$. The entry $M_{i,j}$ is 1 if $x_i \le x_j$, and 0 otherwise.

Consider two posets $(X,\le)$ and $(Y,\le)$. The Cartesian product $X \times Y$ is the poset with the relation
$\le$ defined by
 $(a,b) \le (c,d)$ if and only if $a \le c$ and $b \le d$.

Let us recall that function $f: X \to Y$ is monotone if for any elements $x,y \in X$ such that $x \leq y$, we have $f(x) \leq f(y)$.

We denote the set of all monotone functions from $X$ to $Y$ by $Y^X$. A partial order on $Y^X$ can be defined as follows: given two functions $f,g \in Y^X$, we say that $f \leq g$ if and only if $f(x) \leq g(x)$ for all $x \in X$.

In our algorithms for counting fixed points in $D_n$ under all $\pi \in S_n$, we utilize the following lemmas:

\begin{lemma}\label{L1}\cite{szepietowski} For three posets $R,S,T$:

(1) If $S$ and $T$ are disjoint, then the poset $R^{S+T}$ is isomorphic to $R^S\times R^T$.

(2) The poset $R^{S\times T}$ is isomorphic to $(R^S)^T$ and to $(R^T)^S$.
\end{lemma}

\begin{lemma}\label{L2} \cite{szepietowski}
\leavevmode
\begin{itemize}
\item[(a)] \quad $B^{k+m}=B^{k}\times B^{m}$
\item[(b)] \quad $D_{k+m}=(D_k)^{B^m}$
\end{itemize}
\end{lemma}

These lemmas, in various formulations, are standard tools in the literature for counting monotone Boolean functions. For example, Wiedemann \cite{Wied1991} utilized the isomorphism $D_8 = (D_6)^{B^2}$ to calculate the eighth Dedekind number.

\section{Methodology and results}

To efficiently calculate $\phi_9(\pi)$ for all cases, we employ several algorithms, the outlines of which are described in the following section. For a more detailed explanation of the methodology and examples of these algorithms, please refer to \cite{cz, pawelski, szepietowski}.

\subsection{Basic algorithm}\label{basic}
Let $B^n(\pi)$ denote the poset of orbits of $B^n$ under $\pi \in S_n$. 
Two orbits $C_1$ and $C_2$ are in the relation $C_1 \le C_2$ if and only if $c_1\le c_c$ for
some $c_1\in C_1$ and $c_2\in C_2$.

\begin{center}
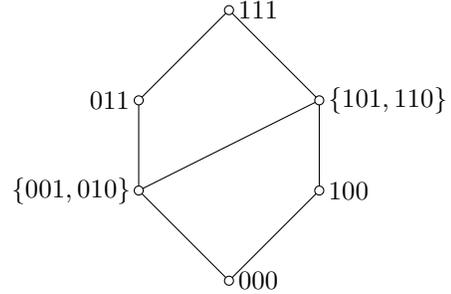

\begin{tikzpicture}[scale=0.6]
\draw (0,0) -- (-2,2);
\draw (-2,2) -- (-2,4);
\draw (0,0) -- (2,2);
\draw (2,2) -- (2,4);
\draw (-2,4) -- (0,6);
\draw (2,4) -- (0,6);
\draw (-2,2) -- (2,4);
\draw[black, fill=white] (0,0) circle (.1cm) node[right] {$000$};
\draw[black, fill=white] (-2,2) circle (.1cm) node[left] {$\{001,010\}$};
\draw[black, fill=white] (2,2) circle (.1cm) node[right] {$100$};
\draw[black, fill=white] (-2,4) circle (.1cm) node[left] {$011$};
\draw[black, fill=white] (2,4) circle (.1cm) node[right] {$\{101,110\}$};
\draw[black, fill=white] (0,6) circle (.1cm) node[right] {$111$};
\end{tikzpicture}
\captionof{figure}{$B^3((12))$ - poset of orbits in $B^3$ under $\pi = (12)$.}
\end{center}

In \cite[Section 3.1]{pawelski}, we show that $\Phi_n(\pi) = B^{B^n(\pi)}$. In other words, for each function in $\Phi_n(\pi)$, bits in $B^n$ are set to the same value on each orbit in $B^n$ under $\pi$. We use this isomorphism to calculate $\phi_n(\pi)$ by generating the set of all downsets of $B^n(\pi)$, using the fact that $|B^{B^n(\pi)}| = \phi_n(\pi)$. Unfortunately, since we store all elements of this set, rather than just counting them, the usability of this algorithm is strictly limited by the amount of main memory. The largest number we can compute using this approach is $\phi_9((123)(456789))= 218542866$; thus, in practice, all sets of fixed points with cardinalities smaller than about two hundred million are within reach using this basic method.

To count fixed points with permutations containing at least one disjoint 1--cycle, we can take advantage of the isomorphism: 
$D_{n+k} = (B^{B^n})^{B^k} = (D_n)^{B^k}$ (see Lemmas \ref{L1} and \ref{L2}), and extend it to $\Phi(\pi)$: $\Phi_{n+k}(\pi) = (B^{B^n(\pi)})^{B^k} = \Phi_n(\pi)^{B^k}$.
 
 In \cite[Section 3.2]{pawelski}, we outline the simplest version of this algorithm for the case when $k=1$. Carić and Živković provide a detailed description of this algorithm for the case when $k=2$ \cite[Theorem 3.1]{cz}. In our research, we utilize this isomorphism for up to $n=4$.

In the calculation of $r_8$, there is only one case where the algorithms described above did not work: specifically, $\phi_8((12)(34)(56)(78))$. In the calculation of $r_9$, we have three cases: $\phi_9((12)(34)(56)(78))$, $\phi_9((12)(34)(56)(789))$, and $\phi_9((123)(456)(789))$.

\subsection{Calculation of \texorpdfstring{$\phi_9((12)(34)(56)(789))$}{phi_9((12)(34)(56)(789))}}\label{p1}

To calculate $\phi_9((12)(34)(56)(789))$, we use a method described by Szepietowski \cite[Section 5]{szepietowski}. We will summarize it briefly. Let $P_n$ denote the chain $\{p_1 \leq p_2 \leq \ldots \leq p_n\}$. The poset $B^2((12))$ is isomorphic to the chain $P_3$ and the poset $B^3((123))$ is isomorphic to the chain $P_4$. Thus, $(B^{n})^{P_k}$ is equal to the sum of elements of the $(k-1)$--power of the incidence matrix of $B^{n}$ (also see: \cite[Proposition 4.8]{aigner}). We therefore have $B^{n+k}(\epsilon) = (B^n(\tau))^{B^k(\gamma)} = (B^n(\tau))^{P_{k+1}}$.
We use this property to calculate $\phi_9((12)(34)(56)(789)) = 807900672006$, which is equal to the sum of all elements of the third power of the incidence matrix of $\Phi_6((12)(34)(56))$.
Another example is the calculation of $\phi_9((12)(345)(6789)) = 22062570$, which is equal to the sum of all elements of the third power of the incidence matrix of $\Phi_6((12)(3456))$.

\subsection{Calculation of \texorpdfstring{$\phi_9((12)(34)(56)(78))$}{phi_9((12)(34)(56)(78)}}\label{p2}

This method is detailed in \cite[Algorithm 3]{pawelski} (also see \cite{cz, szepietowski}) as a response to the need for calculating $\phi_8(\pi)$ under $\pi$ with all cycle lengths of 2. Due to Lemma \ref{L2}, we only need to raise four components to one level. Let $\sigma$ represent the permutation (12).

Thus, each function $x \in \Phi_9((12)(34)(56)(78))$ can be combined from those four functions:

\begin{enumerate}
    \item $a \in \Phi_7(12)(34)(56)$;
    \item $b \in D_7, b \geq a$;
    \item $c = \sigma(b), c \geq a$;
    \item $d \in \Phi_7(12)(34)(56), d \geq b \cup c$.
\end{enumerate}

$\Phi_7(12)(34)(56)$ has 12015832 elements, so we can store it in the main memory. The most challenging aspect is storing $D_7$. We accomplish this by storing $R_7$, which contains 490013148 elements that barely fit in the main memory. Then, for each $x \in R_7$, we perform an "unpacking" operation, making use of the $7! = 5040$ possible permutations of input variables. Finally, we obtained $\phi_9((12)(34)(56)(78)) = 17143334331688770356814$.

\subsection{Calculation of \texorpdfstring{$\phi_9((123)(456)(789)$}{9((123)(456)(789))}}\label{p3}

Let $\omega$ represent the permutation (123)(456). Similarly at in the previous case, we can show that each function $x \in \Phi_9((123)(456)(789))$ can be combined from those eight functions:

\begin{enumerate}
    \item $a \in \Phi_6(123)(456)$;
    \item $b \in D_6, b \geq a$;
    \item $c = \omega(b), c \geq a$;
    \item $d = \omega(c), d \geq a$;
    \item $e \in D_6, e \geq b \cup c$;
    \item $f = \omega(e), f \geq b \cup d$;
    \item $g = \omega(f), g \geq c \cup d$;
    \item $h \in \Phi_6(123)(456), h \geq e \cup f \cup g$.   
\end{enumerate}

$D_6$ has 7828354 elements, $\Phi_6(123)(456)$ has 562 elements, so they fit in the main memory. We completed a calculation in about six hours with a result $\phi_9((123)(456)(789) = 221557843276152$.

\subsection{Implementation and results}

The algorithms were implemented in Java and Rust. The basic algorithm (Section \ref{basic}) was implemented in Java, while the cases from Sections \ref{p1}, \ref{p2}, and \ref{p3} were implemented in Rust. We used a machine with 32 Xeon threads. The total computation time is challenging to estimate accurately, as we calculated each of the 29 cases separately. However, with all algorithms implemented, we estimate the total computation time to be about 10 days. The simplest cases are computed instantly, while the most challenging ones (for example, $\phi_9(12)(34)$) take up to 2 days.

\begin{table}[H]
\centering
    {\renewcommand{\arraystretch}{1.5}
    \begin{tabular}{c*{3}{c}c}
    $\pi_i$ & $\mu_i$ & $\phi_9(\pi_i)$ \\
    \hline
        (12)              & 36      & 16278282012194909428324143293364  \\
        (123)             & 168     & 868329572680304346696  \\
        (1234)            & 756     & 5293103318608452  \\
        (12345)           & 3024    & 26258306096  \\
        (123456)          & 10080   & 2279384919  \\
        (1234567)         & 25920   & 3268698  \\
        (12345678)        & 45360   & 1144094  \\
        (123456789)       & 40320   & 97830  \\
        (12)(34)          & 378     & 107622766375525877620879430  \\
        (12)(345)         & 2520    & 5166662396125146  \\
        (12)(3456)        & 7560    & 323787762940974  \\
        (12)(34567)       & 18144   & 70165054  \\
        (12)(345678)      & 30240   & 547120947  \\
        (12)(3456789)     & 25920   & 80720  \\
        (123)(456)        & 3360    & 7107360458115201  \\
        (123)(4567)       & 15120   & 92605092  \\
        (123)(45678)      & 24192   & 197576  \\
        (123)(456789)     & 20160   & 218542866  \\
        (123)(456)(789)   & 2240    & 221557843276152  \\
        (1234)(5678)      & 11340   & 503500313130  \\
        (1234)(56789)     & 18144   & 10182  \\
        (12)(34)(56)      & 1260    & 328719964864138799170044  \\
        (12)(34)(567)     & 7560    & 14037774553676  \\
        (12)(34)(5678)    & 11340   & 66031909836340 \\
        (12)(34)(56789)   & 9072    & 3710840  \\
        (12)(345)(678)    & 10080   & 866494196253 \\
        (12)(345)(6789)   & 15120   & 22062570  \\
        (12)(34)(56)(78)  & 945     & 17143334331688770356814  \\
        (12)(34)(56)(789) & 2520    & 807900672006  \\
        \hline
    \end{tabular}}
    \[\sum_{i=2}^{k} \mu_i \phi_9(\pi_i) = 586059264378237446637837193706034.\]
    \caption{Values of all $\phi_9(\pi)$ (excluding the identity).} \label{tab:r9results}
    \end{table}

\section{Calculation of \texorpdfstring{$r_9$}{r(9)}}

In April 2023, two research teams independently reported the following value of $d_9$ \cite{jakel, hirtum}:

$$ d_9 = 286386577668298411128469151667598498812366.$$
We can now finally make direct use of Equation \ref{burnside}, obtaining the following value: 

$$ r_9 = 789204635842035040527740846300252680.$$

\appendices

\section{\texorpdfstring{$r_n$}{r(n)} -- calculation tables}
In this section, we present the number of fixed points for all cycle types in $D_n$ for $n$ up to 8 and the application of Equation \ref{burnside} using the obtained values. We have recalculated all values in these tables.

\begin{table}[H]
\centering
    {\renewcommand{\arraystretch}{1.5}
    \begin{tabular}{c*{3}{c}c}
    $i$ & $\pi_i$ & $\mu_i$ & $\phi_9(\pi_i)$ \\
    \hline
        1 & (1)              & 1      & 6  \\
        2 & (12)             & 1      & 4  \\
        \hline
    \end{tabular}}
    \[r_2 = \frac{1}{2!} \cdot \sum_{i=1}^{k} \mu_i \phi_2(\pi_i) = 5.\]
    \caption{Calculation of $r_2$.}
    \end{table}

\begin{table}[H]
\centering
    {\renewcommand{\arraystretch}{1.5}
    \begin{tabular}{c*{3}{c}c}
    $i$ & $\pi_i$ & $\mu_i$ & $\phi_9(\pi_i)$ \\
    \hline
        1 & (1)              & 1      & 20  \\
        2 & (12)             & 3      & 10  \\
        3 & (123)            & 2      & 5  \\
        \hline
    \end{tabular}}
    \[r_3 = \frac{1}{3!} \cdot \sum_{i=1}^{k} \mu_i \phi_3(\pi_i) = 10.\]
    \caption{Calculation of $r_3$.}
    \end{table}

\begin{table}[H]
\centering
    {\renewcommand{\arraystretch}{1.5}
    \begin{tabular}{c*{3}{c}c}
    $i$ & $\pi_i$ & $\mu_i$ & $\phi_9(\pi_i)$ \\
    \hline
      1 &  (1)              & 1      & 168  \\
      2 &  (12)             & 6      & 50  \\
      3 &  (123)            & 8      & 15  \\
      4 &  (1234)           & 6      & 8  \\
      5 &  (12)(34)         & 3      & 28  \\
        \hline
    \end{tabular}}
    \[\frac{1}{4!} \cdot \sum_{i=1}^{k} \mu_i \phi_4(\pi_i) = 30.\]
    \caption{Calculation of $r_4$.}
    \end{table}

\begin{table}[H]
\centering
    {\renewcommand{\arraystretch}{1.45}
    \begin{tabular}{c*{3}{c}c}
    $i$ & $\pi_i$ & $\mu_i$ & $\phi_9(\pi_i)$ \\
    \hline
       1 & (1)              & 1       & 7581  \\
       2 & (12)             & 10      & 887  \\
       3 & (123)            & 20      & 105  \\
       4 & (1234)           & 30      & 35 \\
       5 & (12345)          & 15      & 309  \\
       6 & (12)(34)         & 24      & 11  \\
       7 & (12)(345)        & 20      & 35  \\
        \hline
    \end{tabular}}
    \[\frac{1}{5!} \cdot \sum_{i=1}^{k} \mu_i \phi_5(\pi_i) = 210.\]
    \caption{Calculation of $r_5$.}
    \end{table}

\begin{table}[H]
\centering
    {\renewcommand{\arraystretch}{1.45}
    \begin{tabular}{c*{3}{c}c}
    $i$ & $\pi_i$ & $\mu_i$ & $\phi_9(\pi_i)$ \\
    \hline
       1 &  (1)              & 1      & 7828354  \\
       2 &  (12)             & 15     & 160948  \\
       3 &  (123)            & 40     & 3490  \\
       4 &  (1234)           & 90     & 494 \\
       5 &  (12345)          & 144    & 64  \\
       6 &  (123456)         & 120    & 44 \\
       7 &  (12)(34)         & 45     & 24302  \\
       8 &  (12)(345)        & 120    & 490  \\
       9 &  (12)(3456)       & 90     & 324  \\
       10 & (123)(456)       & 40     & 562 \\
       11 & (12)(34)(56)     & 15     & 8600  \\
        \hline
    \end{tabular}}
    \[r_6 = \frac{1}{6!} \cdot \sum_{i=1}^{k} \mu_i \phi_6(\pi_i) = 16353.\]
    \caption{Calculation of $r_6$.}
    \end{table}

\begin{table}[H]
\centering
    {\renewcommand{\arraystretch}{1.45}
    \begin{tabular}{c*{3}{c}c}
    $i$              & $\pi_i$ & $\mu_i$ & $\phi_7(\pi_i)$ \\
    \hline
        1    & (1)              & 1     & 2414682040998  \\
        2    & (12)             & 21    & 2208001624  \\
        3    & (123)            & 70   & 2068224  \\
        4    & (1234)           & 210   & 60312  \\
        5    & (12345)          & 504  & 1548  \\
        6    & (123456)         & 840  & 766  \\
        7    & (1234567)        & 720  & 101  \\
        8    & (12)(34)         & 105   & 67922470  \\
        9    & (12)(345)        & 420  & 59542  \\
        10   & (12)(3456)       & 630  & 26878  \\
        11   & (12)(34567)      & 504  & 264  \\
        12   & (123)(456)       & 280  & 69264  \\
        13   & (123)(4567)      & 420  & 294  \\
        14   & (12)(34)(56)     & 105   & 12015832  \\
        15   & (12)(34)(567)    & 210  & 10192  \\
        \hline
    \end{tabular}}
    \[r_7 = \frac{1}{7!} \cdot \sum_{i=1}^{k} \mu_i \phi_7(\pi_i) = 490013148.\]
    \caption{Calculation of $r_7$.\cite{pawelski}}
    \end{table}

\begin{table}[H]
\centering
    {\renewcommand{\arraystretch}{1.5}
    \begin{tabular}{c*{3}{c}c}
        $i$              & $\pi_i$ & $\mu_i$ & $\phi_8(\pi_i)$ \\
    \hline
        1    & (1)              & 1     & 56130437228687557907788  \\
        2    & (12)             & 28    & 101627867809333596  \\
        3    & (123)            & 112   & 262808891710  \\
        4    & (1234)           & 420   & 424234996  \\
        5    & (12345)          & 1344  & 531708  \\
        6    & (123456)         & 3360  & 144320  \\
        7    & (1234567)        & 5760  & 3858  \\
        8    & (12345678)       & 5040  & 2364  \\
        9    & (12)(34)         & 210   & 182755441509724  \\
        10   & (12)(345)        & 1120  & 401622018  \\
        11   & (12)(3456)       & 2520  & 93994196  \\
        12   & (12)(34567)      & 4032  & 21216  \\
        13   & (12)(345678)     & 3360  & 70096  \\
        14   & (123)(456)       & 1120  & 535426780  \\
        15   & (123)(4567)      & 3360  & 25168  \\
        16   & (123)(45678)     & 2688  & 870  \\
        17   & (1234)(5678)     & 1260  & 3211276  \\
        18   & (12)(34)(56)     & 420   & 7377670895900  \\
        19   & (12)(34)(567)    & 1680  & 16380370  \\
        20   & (12)(34)(5678)   & 1260  & 37834164  \\
        21   & (12)(345)(678)   & 1120  & 3607596  \\
        22   & (12)(34)(56)(78) & 105   & 2038188253420  \\
        \hline
    \end{tabular}}
    \[r_8 = \frac{1}{8!} \cdot \sum_{i=1}^{k} \mu_i \phi_6(\pi_i) = 1392195548889993358.\]
    \caption{Calculation of $r_8$.\cite{pawelski}}
    \end{table}

\vspace{1mm}

\appendices

\begin{IEEEbiographynophoto}{Bart\l{}omiej Pawelski} works at University of Gdansk, Poland.
\end{IEEEbiographynophoto}

\end{document}